\documentclass[12pt]{article}
\usepackage{amsmath,amsfonts,latexsym}

\newcommand{\g}{\mbox{$\mathfrak g$}}
\newcommand{\p}{\mbox{$\mathfrak p$}}
\begin{document}
\begin{center}
{\bf Strongly Semistable Bundles on a Curve over a finite field} \\
{\bf S. Subramanian}\\[5mm]
\end{center}
\section*{Abstract}
	We show that a principal G bundle on a smooth projective curve
over a finite field is strongly semistable if and only if it is defined
by a representation of the fundamental group scheme of the curve into G.
\\[10mm]
MSC 2000 classification:14-XX
\\[5mm]
Keywords:semistable bundles,fundamental group scheme.
\\[5mm]
\section*{1. Introduction}

\quad Let $G$ be a semisimple algebraic group and $X$ a smooth projective 
curve,  both defined over a finite field $k$ of characteristic $p$.  Let $E 
\to  X$ be a principal $G$-bundle on $X$.  $E$ is said to be  stable
(resp. 
semistable) if for every reduction of structure group $E_P \subset E$ to a 
parabolic subgroup $P$ of $G$, we have
$${\rm deg} \ E_P (\p) < ({\rm resp.} \le) \ {\rm deg} \ E_G 
(\g ) $$
where $\g$  and $\p$ denote the lie algebras of $G$ and $P$
respectively,and $E_P(\p)$(respectively $E_G(\g)$)denotes the Lie algebra
bundle associated to $E_P$(respectively $E_G$) by the adjoint action of
$P$(respectively $G$)on $\p$(respectively $\g$). 
When the base field is the field of complex numbers, the theorem  of
Narasimhan-Seshadri-Ramanathan  says that a stable $G$-bundle is defined by a
unitary representation of the topological fundamental group of $X$.  In this 
article, we show that in the case when $k$ is a finite field, a strongly 
semistable principal $G$-bundle on $X$ (recall that a bundle is said to be 
strongly semistable if all of its Frobenius pull backs are semistable) is 
defined by a representation of the fundamental group scheme $\pi(X) \to 
G$.  Here, by fundamental group scheme, we mean the proalgebraic group scheme
introduced by M. Nori in [3], [4].
	The author wishes to thank the referee for his helpful comments on
an earlier version of this paper. 
\section*{2.Preliminaries}
\quad Let $k$ be a finite field of characteristic $p$ and let $X$ be a smooth
projective curve defined over $k$.  We assume that $X$ has a rational point 
(say $x$) defined over $k$.  A vector bundle $V \to X$ (over $k$) of degree 
zero is said to be essentially finite if there exists a finite group
scheme 
$H$ over $k$ and a principal $H$-bundle $\pi : E_0 \to X$ on $X$, such that
$\pi^*V$ is trivial on $E_0$ (see Proposition(3.10),page 83,in [5]).  Let
${\cal C}$ denote the 
category of all essentially finite vector bundles defined on $X$ over $k$.  We 
define a fibre functor $T$ on ${\cal C}$ by defining $T(V) = V_x$, the fibre 
of $V$ at $x$, for every $V \in {\cal C}$.  With this fibre functor,
${\cal C}$  is a Tannaka category over $k$ (see section 2.3 in [5]).  The
associated
affine
group scheme $\pi(X, x)$ is called the fundamental group scheme of $X$ 
over $k$. Further, there is a principal $\pi (X, x)$ bundle ${\cal X} \to
X$ over 
$k$ (see section 2.3,page 84,in [5]).  Given a semisimple algebraic group
$G$ over $k$, and an 
algebraic group homomorphism $\pi(X, x) \to G$ (such a homomorphism factors 
through a finite group scheme quotient of $\pi (X, x))$, there is associated
a principal $G$-bundle $E\to X$, associated to ${\cal X} \to X$ by the
homomorphism 
$\pi(X, x) \to G$. 

We now consider a principal $G$ bundle $E \to X$, which we assume to be 
strongly semistable.   Since $G$  is assumed to be semisimple, we can find a 
$G$-module $V$ such that $G \subset SL(V)$.  The associated vector bundle
$E(V)$ is strongly semistable by the theorem of Ramanan-Ramanathan  
(see Theorem 3.23 in [6] ).   We form the category ${\cal M}$ whose
objects are subquotients of degree zero of finite direct sums 
$${\displaystyle{\bigoplus_{i}}} E(V)^{\oplus n_i} \oplus 
E(V) ^{* \oplus m_i}$$
where $n_i, m_i$ are nonnegative integers.The homomorphisms in ${\cal M}$
are all homomorphisms of the respective vector bundles on $X$.
Given a vector bundle $W$ in 
${\cal M}$, we can associate to it its fibre $W_x$ at $x$, and with this
as 
the fibre functor, ${\cal M}$ is a Tannaka category over $k$ (see Section 
4 in [1]).  
There is an associated affine group scheme $M$, called the monodromy group 
scheme of $E$.Further,if $Vect(X)$ denotes the category of vector bundles
on $X$,then there is a naturally defined functor
$\phi:{\cal M} \rightarrow Vect(X)$ which associates to an object
of ${\cal M}$ the vector bundle on $X$ that it represents.The functor
$\phi$ satisfies the axioms of Lemma(2.1),page 77,in [5].Hence,by
Proposition(2.5),page 78,in [5],we obtain a principal $M$-bundle
$E_M \to X$ on $X$, defined over $k$ (see [1],section 4).  Further, since
$k$ is a finite field, $H^1(k,G)$ is trivial (see Theorem 1',section
2.3,Chapter 3,page 132 in [7]), and hence we can assume that there is a
$k$ rational point
in the fibre $E_x$ of $E$ at $x$.  With this choice of point in $E_x$, we 
obtain an embedding $M \subset G$, and a reduction of structure group 
$E_M \subset E$ (see section 4 in [1] for details).  The monodromy group
scheme $M$ and 
the bundle $E_M$ do not depend on the choice of representation $G \subset 
SL(V)$,and the bundle $E_M$ is characterised as the minimal degree zero
reduction of structure group of $E$(see Proposition 4.9,Corollary
4.10 and Proposition 4.11 in [1]).  
\section*{3.Main Theorem}  Let $k$ be a finite field and let $X$ be a
smooth projective curve over $k$.  We assume that $X$ has a rational point
$x$ over $k$.  Let 
$G$ be a semisimple algebraic group over $k$.  Let $F$ denote the Frobenius 
morphism on $X$. 
\paragraph*{Definition \ (3.1)}  A  principal $G$ bundle $E \to X$ on $X$
over $k$ is 
said to be semistable if, for every reduction of structure group $E_P 
\subset
E$ to a parabolic subgroup $P$ of $G$, we have
$$ {\rm degree \ } E_P(\p) \le 0.$$ 
The principal bundle is said to be strongly semistable if $F^{n^*}E$ is 
semistable for all $ n \ge 1$.  We can now prove 
\paragraph*{Theorem (3.2)} Let $E\to X$ be a strongly semistable $G$-bundle 
on $X$ over $k$.  Then $E\to X$ is defined by a homomorphism $\pi(X,x) \to
G$. 
This homomorphism $\pi(X,x) \to G$ is unique upto inner conjugation
by an element of $G(k)$.  Further, if $G'$ is another 
semisimple algebraic group over $k$ and $G \to G'$ a homomorphism of
algebraic groups, then the strongly semistable $G'$ 
bundle $E_{G'} \to X$ associated to $E \to X$ is defined by the composite 
homomorphism $\pi(X,x) \to G \to G'$.  
\paragraph*{Proof.}  By hypothesis, the bundles $F^{n*}E, n \ge 1$ are all 
semistable.  Since $k$ is a finite field, by boundedness, the set $\{ F^{n*} E,
n \ge 1\}$ is a finite set.  Hence $F^{n*} E \cong F^{(n+ \ell)*}E$ for 
some positive integers $n$ and $\ell$.  We now consider the monodromy group
scheme $M_n$, and the monodromy bundle $E_{M_n}$, of the bundle $F^{n*} E$ for 
$n \ge 0$ (with the convention that $F^{0 *} E = E)$.  Let $F_n$ denote
the 
Frobenius morphism of the group scheme $M_n$.  The isomorphism  $F^{n*}E
\cong
F^{(n+\ell)*} E$ implies that $F^{\ell}_n(M_n)\cong M_{n+\ell} \cong M_n$,
and 
further the monodromy bundles $E_{M_n}$ and $E_{M_{n+\ell}}$ are isomorphic.  
We now observe that the isomorphism 
$$F_n^{\ell}(M_n) \cong M_{n+\ell}\cong M_n$$ 
implies that the group scheme $M_n$ is finite and reduced.Then 
by Lemma (4.12) in [1], it follows that $M_{0,red}$, the reduced part of 
the group scheme $M_0$(where $M_0$ denotes the monodromy group scheme of
the bundle $E$), is isomorphic to $M_n$.Also,if $F_0$ denotes the
Frobenius morphism of $M_0$,we obtain $F^{n}_0(M_0) \cong M_n$.Hence $M_0$
is a  finite group scheme over $k$.  The monodromy bundle $E_{M_0}$ is
a principal 
$M_0$ bundle on $X$, and the reduction of structure group $E_{M_0}\subset E$ 
implies that the pull back of the bundle $E$ to $E_{M_0}$ is trivial.  The 
principal $M_0$ bundle $E_{M_0} \to X$ defines a quotient $\pi(X, x) \to 
M_0$ (see [5]), and the composite homomorphism 
$$\pi (X, x)  \to M _0 \to G$$ 
is the homomorphism we want defining the principal $G$-bundle $E$.  

For the second part, let $G \to G'$ be a homomorphism of semisimple groups  
over $k$.  Since $G$
is semisimple, there are three cases to consider:
\begin{enumerate}
\item[Case (i)] $G \subset G'$ is a closed immersion: in this case, a
faithful
representation of $G'$ is also a faithful representation of $G$, and the 
result follows from the fact that the monodromy group scheme and the monodromy
bundle of $E$ are independent of the choice of representation. 
\item[Case (ii)] $G \to G'$ is a central isogeny: we consider the case 
$G' = G/Z$, where $Z$ is the centre of $G$.  Let $M$ be the monodromy group 
scheme  of $E$ and $E_M$ the monodromy bundle.  Let $M'$ be the image of 
$M$ in $G'$ and $E_{M'}$ the bundle associated to $E_M$ by the projection 
$M \to M'$.  Then $E'$ has a reduction of structure group $E_{M'} \subset
E'$ to $M'$, of degree zero(see Definition 4.7 in [1]).  If $H$ is a group
scheme over $k$ and $E_H
\subset E'$ a reduction of structure group of $E'$ of degree zero to $H$, 
with $H \subset M'$, then taking $\widetilde{H}$ to be the inverse image 
of $H$ in $G$, we obtain a principal $\widetilde{H}$ bundle $E_{\widetilde{H}}$
and a reduction of structure group $E_{\widetilde{H}} \subset E$ of degree 
zero  with $\widetilde{H} \subset M$.  
This contradicts the minimality of $M$ and $E_M$ (see Corollary (4.10) in 
[1]).  It  follows that $M'$ is the minimal subgroup scheme of $G'$ with 
a degree zero reduction of structure group to $M'$.  Hence, again by 
Corollary (4.10) in [1], it follows that $M'$ is the monodromy group
scheme
of $E'$ and $E_{M'}$ is the monodromy bundle of $E'$. 
\item[Case(iii)]$G \to G'$ is purely inseparable:to treat this case,it is
enough to consider the Frobenius morphism $F:G \to G$,as the general
case is what is usually called a Frobenius sandwich.Given the principal
G bundle $E \to X$,the bundle induced by the Frobenius $F:G \to G$ is
the pullback $F^*E$(where the pullback is taken under the Frobenius
morhism of $X$).Let $M$ be the monodromy groupscheme of $E$ and
$E_M \subset E$ the monodromy reduction.We know that $M$ is a finite
groupscheme.Let $F_M$ denote the Frobenius morphism of $M$ and let
$M'$ be the image of $M$ under $F_M$.The $M$ bundle $E_M$ induces
a $M'$ bundle $E_{M'}$ and a reduction of structure group
$E_{M'} \subset F^*E$.If there exists a proper subgroupsheme 
$N \subset M'$ and a reduction of structure group
$E_N \subset E_{M'}$, then the inverse image of $E_N$ under the
projection $E_M \to E_{M'}$ would define a reduction of structure
group of $E_M$,contradicting the minimality of $E_M$.Therefore
$E_{M'}$ is a minimal reduction of structure group of $F^*E$,and
hence is the monodrmy reduction. 
\end{enumerate}
It now follows that the monodromy of the bundle $E'$ is 
given by the composite homomorphism $\pi(X, x) \to G \to G'$.  \hfill{Q.E.D.} 

\section*{4.When $X$ does not have a rational point}
   Let $X$ be a smooth projective curve over a finite field $k$ as before
and let $E \to X$ be a principal G bundle on $X$ which is strongly
semistable.Let $l$ be a finite Galois extension of $k$ over which $X$
acquires a rational point.Let $X_l$ denote the base change of $X$ to
$Spec(l)$ and $E_l$ the base change of $E$ to $Spec(l)$.Then by $section
3$
above,$E_l$ is defined by a homomorphism $\pi(X_l) \to G_l$,where 
$\pi(X_l)$ denotes the fundamental group sheme of $X_l$ and $G_l$
denotes the base change of $G$ to $Spec(l)$.
    We have the following exact sequence of group schemes over $Spec(k)$,
$$1 \to \pi(X_l,x) \to \pi(X,x) \to Gal(l/k) \to 1$$
where $x$ is a closed point of $X$ with residue field $l$,so it is
a rational point over $l$,$\pi$ is the fundamental group scheme,and
$Gal(l/k)$ is the Galois group of $l$ over $k$.Here $\pi(X_l,x)$,which
is a group scheme over $Spec(l)$,is regarded as a group scheme over
$Spec(k)$ by the projection $Spec(l) \to Spec(k)$.The homomorphism
$\pi(X_l,x) \to G_l$ factors through a finite group scheme quotient,
say $H$,of $\pi(X_l,x)$.Since $E_l$ is actually defined over $k$,it
follows that the monodromy group scheme $M$ of $E_l$,and the monodromy
reduction $E_M \subset E_l$, remain invariant under the action of
$Gal(l/k)$ on $E_l$.Hence the subgroupscheme $H \subset G_l$ remains
invariant under the action of $Gal(l/k)$.Every groupsheme quotient
of $\pi(X_l,x)$ invariant under the action of $Gal(l/k)$ on $\pi(X_l,x)$
descends to a quotient of $\pi(X,x)$.Hence $H$ is a quotient of
$\pi(X,x)$ and we obtain the homomorphism $\pi(X,x) \to G$ defining
$E$ on $X$ over $k$. 	
\section*{References}
\begin{enumerate}
\item Indranil Biswas, A.J. Parameswaran and S. Subramanian: Monodromy
group 
for a strongly semistable principal bundle over a curve, {\it Duke Math. 
J.} {\bf 132}, No.1, 1-48 (2006).
\item H.Esnault,P.H.Hai,and X.Sun,On Nori's fundamental group scheme,
arXiv:math.AG/0605645 
\item H. Lange and U. Stuhler: Vectorbundel auf Kurven und Darstellungen 
der algebraischen Fundamentalgruppe, {\it Math. Z.} {\bf 156}, No.1, 73-83 
(1977). 
\item M.V. Nori: On the representations of the fundamental group scheme, 
{\it Comp. Math.}, {\bf 33}, 29-41 (1976). 
\item M.V. Nori: The fundamental group scheme, {\it Proc. Ind. Acad. Sci. 
(Math. Sci.)}, {\bf 91}, 73-122 (1982).
\item S.Ramanan and A.Ramanathan:Some remarks on the instability flag,{\it
Tohoku Math.J.}(2),{\bf 36},269-291(1984).
\item J.P. Serre: Galois cohomology,English translation by
P.Ion,Springer,1997\\[10mm]
\end{enumerate}
School of Mathematics \\
Tata Institute of Fundamental Research \\
Homi Bhabha Road, Colaba \\
Mumbai 400 005 \\
India\\
E-mail:subramnn@math.tifr.res.in   
\end{document}